\documentclass[10pt]{article} 
\usepackage[utf8]{inputenc}
\usepackage[english]{babel}
\usepackage[]{amsthm} 
\usepackage[]{amssymb} 
\usepackage[]{amsmath}
\usepackage{faktor}
\usepackage{xcolor}
\usepackage{dsfont}
\usepackage[shortlabels]{enumitem}
\usepackage{tikz-cd}

\usepackage{graphicx}
\graphicspath{ {./} }

\usepackage{csquotes}

\usepackage{geometry}
\newcommand{\BC}{\mathbb{C}}
\newcommand{\BR}{\mathbb{R}}
\newcommand{\ot}{\otimes}

\newcommand{\HH}{\mathcal{H}}
\newcommand{\PP}{\mathcal{P}}

\newtheorem{thm}{Theorem}

\newtheorem{prop}[thm]{Proposition}

\theoremstyle{definition}
\newtheorem{defn}[thm]{Definition}
\newtheorem{rem}[thm]{Remark}
\newtheorem{ex}[thm]{Example}
\newtheorem{prop normal}[thm]{Proposition}
\newtheorem{conj}[thm]{Conjuncture}

\usepackage{alphalph}

\title{Families of Perfect Tensors}
\author{Runshi Geng}

\begin{document}
\maketitle
\begin{abstract}
\noindent Perfect tensors are the tensors corresponding to the absolutely maximally entangled states, a special type of quantum states of interest in quantum information theory. We establish a method to compute parameterized families of perfect tensors in $(\BC^d)^{\ot 4}$ using exponential maps from Lie theory. With this method, we find explicit examples of non-classical perfect tensors in $(\BC^3)^{\ot 4}$. In particular, we answer an open question posted by {\.Z}yczkowski et al.
\end{abstract}

\section{Introduction}
The entanglement of quantum states plays a central role in quantum information theory. For bipartite states, the entanglement is usually measured by the von Neumann entropies of the partial traces of the states \cite{MR1418618}, but its generalization for multipartite states is more complicated. One is interested in an extremal class of multipartite states: the absolutely maximally entangled (AME) states -- the pure states that are maximally entangled for every bipartitions of the system. The AME states are closely related to quantum secret sharing \cite{ame}, quantum error-correcting \cite{MR3370186} and quantum repeaters \cite{alsina2021absolutely}.

The existence of AME states in $(\BC^d)^{\ot n}$, denoted as AME$(n,d)$, is an open-problem \cite{openproblem}. For $n=2$ and $3$, the Bell and GHZ states are AME states. For $n=4$, it was proved that AME$(4,2)$ does not exist \cite{clarisse2005entangling,colbourne2007handbook}, and AME$(4,d)$ can be obtained by orthogonal Latin squares (also known as Graeco-Latin squares) of order $d$, denoted as OLS($d$) \cite{colbourne2007handbook}. OLS($d$)'s have been constructed for all $d>2$ except $d=6$ \cite{OLS,MOLS}, which implies the existence of AME$(4,d)$ for such $d$. Although OLS(6) was proved not existing \cite{OLS6}, an explicit example of AME(4,6) is found in \cite{MR4397007}.

The tensors corresponding to the AME states are perfect tensors \cite{MR3370186}. An open question is asked in \cite{zyczkowski20229}: for $d=3,4$ and $5$, are there nonzero perfect tensors in $(\BC^d)^{\ot 4}$ that are non-classical, i.e., not locally equivalent to a solution corresponding to an OLS($d$)? In this paper, we establish a method with the exponential maps from Lie theory to find parameterized families of perfect tensors in $(\BC^d)^{\ot 4}$ for $d>2$ and $d\neq 6$. Particularly we find explicitly parameterized families of perfect tensors in $(\BC^3)^{\ot 4}$ and we proved that there are 4-dimensional subsets in those families consisting of non-classical perfect tensors.

\section{Preliminaries}
\subsection{The Varieties of Perfect Tensors}
Let $\HH:=\HH_1\ot\cdots\ot\HH_n$, where each $\HH_i$ is a $d$-dimensional complex Hilbert space with the inner product $\langle \Phi|\Psi\rangle:=\Phi^{\dagger}\Psi$. 

Write $[n]:=\{1,2,\cdots,n\}$. A \textbf{bipartition} of $[n]$ is a pair of non-empty subsets $I,J\subset [n]$ such that $I\cup J= [n]$ and $I\cap J=\emptyset$. Identifying $\HH_I$ with its dual space $\HH_I^*$, any tensor $\Phi\in\HH$ can be regarded as a linear map $\Phi_{I,J}:\HH_I\to \HH_J$, where $\HH_I:=\bigotimes_{i\in I} \HH_i$ and $\HH_J:=\bigotimes_{j\in J} \HH_j$.

\begin{defn}[\cite{MR3370186}]
A tensor $\Phi\in \HH$ is called a \textbf{perfect tensor} if for every bipartition $(I,J)$ of $[n]$ with $|I|\leq |J|$, $\Phi_{I,J}$ is proportional to an isometry.
\end{defn}

For even $n$, let $\PP(n,d)$ be the set of tensors $\Phi\in\HH$ such that for every bipartition $(I,J)$ of $[n]$ with $|I|=|J|=\frac{n}{2}$ , $\Phi_{I,J}$ is unitary.

\begin{prop}[\cite{MR3370186}]
When $n$ is even, $\Phi\in\HH$ is perfect if and only if it is proportional to a tensor in $\PP(n,d)$.
\end{prop}

That is, the set of all perfect tensors in $\HH$ is just the cone over the real algebraic variety $\PP(n,d)$, which we will mainly study for the rest of this paper. 

For $n=4$, write $\HH=A\ot B\ot C\ot D:=(\BC^d)^{\otimes 4}$. There are 3 flattening maps:
\begin{align*}
    F_1: \HH\to (A\ot B)\ot(C\ot D),\\
    F_2: \HH\to (A\ot C)\ot(B\ot D),\\
    F_3: \HH\to (A\ot D)\ot(C\ot B).
\end{align*}

Let $U_1, U_2$ and $U_3$ be the unitary groups in $(A\ot B)\ot(C\ot D),(A\ot C)\ot(B\ot D)$ and $(A\ot D)\ot(C\ot B)$ respectively. Then
\begin{equation}\label{3int}
    \PP(4,d)=F_1^{-1}(U_1)\cap F_2^{-1}(U_2)\cap F_3^{-1}(U_3).
\end{equation}

Generally, when $n$ is even, there are $\frac{1}{2} {n\choose n/2}$ flattening maps $F_i$'s and unitary groups $U_i$'s, then $\PP(n,d)=\bigcap_{i} F_i^{-1}(U_i)$.

\subsection{Orthogonal Latin Squares}
\begin{defn}[\cite{colbourne2007handbook}]
An \textbf{orthogonal Latin square} of order $d$ (OLS($d$)) is an ordered pair of two $d\times d$ matrices $A=(a_{ij})_{i,j}$, $B=(b_{ij})_{i,j}$ with entries in $[d]$, such that each integer in $[d]$ appears exactly once in each column and each row of $A$ and $B$, and all $d^2$ pairs $(a_{ij},b_{ij})$ are distinct.
\end{defn}

Any OLS($d$) gives us a tensor $\Phi=\sum\Phi_{ijkl}|ijkl\rangle\in(\BC^d)^{\ot 4}$ via $\Phi_{ijkl}=\delta_{a_{kl},i}\delta_{b_{kl},j}$. Then $F_i(\Phi)$ is a permutation matrix for $i=1,2,3$, so $\Phi\in\PP(4,d)$. On the contrary, any tensor $\Phi\in(\BC^d)^{\ot 4}$ such that $F_i(\Phi)$ is a permutation matrix for $i=1,2,3$ yields an OLS($d$). \cite{MR4397007}

\begin{ex}[\cite{MR4397007}]\label{OLS3} An example of OLS(3), written as a $3\times3$ table whose $(i,j)$-th cell is the pair $(a_{ij},b_{ij})$:
\begin{center}
\begin{tabular}{ |c|c|c| }
\hline
 2,3 & 3,1 & 1,2\rule[-1ex]{0pt}{4ex}\\ 
\hline
 3,2 & 1,3 & 2,1\rule[-1ex]{0pt}{4ex}\\
\hline
 1,1 & 2,2 & 3,3\rule[-1ex]{0pt}{4ex}\\
\hline
\end{tabular}
\end{center}
\end{ex}
corresponds to the tensor $\Phi=|1123\rangle+|1231\rangle+|1312\rangle+|2132\rangle+|2213\rangle+|2321\rangle+|3111\rangle+|3222\rangle+|3333\rangle$.

\subsection{Exponential Maps}
For any matrix Lie group $G$ and any element $g\in G$, let $T_gG$ be the tangent space of $G$ at $g$, and $\mathrm{exp}_g:T_gG\rightarrow G$ be the exponential map at $g$. One has the left translation map $L_g$ and its differential $dL_g$: 
\begin{align*}
    L_g: G&\longrightarrow G, \\
     g'&\longmapsto gg',\\
    dL_g: T_{\mathrm{Id}} G&\longrightarrow T_g G, \\
     X&\longmapsto gX.
\end{align*}

The exponential map at the identity $\mathrm{exp}_{\mathrm{Id}}:T_{\mathrm{Id}} G\to G$ agrees with the matrix exponential, that is:
$$\mathrm{exp}_{\mathrm{Id}}(X)=\sum_{i=0}^{\infty} \frac{X^i}{i!}.$$
We obtain the commutative diagram
$$\begin{tikzcd}
T_{\mathrm{Id}}G \arrow{r}{dL_g} \arrow[swap]{d}{\mathrm{exp}_{\mathrm{Id}}} & T_{g}G \arrow{d}{\mathrm{exp}_g} \\
G \arrow{r}{L_g} & G
\end{tikzcd}
$$
which implies that
\begin{equation}\label{exp}
    \mathrm{exp}_g(X)=g\mathrm{exp}_{\mathrm{Id}}(g^{-1}X).
\end{equation}

\section{Parameterized Families in $\PP(4,d)$}
The equation (\ref{3int}) shows $\PP(4,d)$ is determined by numerous quadratic equations, and it is hard to find a complete solution and understand its geometric structure. However, if there is a sparse tensor $\Phi\in\PP(4,d)$, it is much easier to find the intersection of the tangent spaces $T_{\Phi}F_i^{-1}(U_i)$. Then by applying the exponential maps to tangent vectors in a subset of $T:=\bigcap_{i=1}^3 T_{\Phi}F_i^{-1}(U_i)$, we can obtain some curves of, and possibly higher dimensional manifolds of, tensors in $\PP(4,d)$.

The terminology \lq\lq exponential map" might not make sense on $T$ or $T_{\Phi}\PP(4,d)$ if $\PP(4,d)$ is not smooth at $\Phi$. Let $Z_1,\cdots,Z_m\subset\PP(4,d)$ be all maximal irreducible subvarieties of $\PP(4,d)$ that are smooth at $\Phi$. Then $\bigcup_j T_{\Phi}Z_j$ is the subset of $T_{\Phi}\PP(4,d)$ where the exponential map exists. We have 
$$T \supset T_{\Phi}\PP(4,d) \supset \bigcup_{j} T_{\Phi}Z_j.$$

For $i=1,2,3$, let $\mathrm{exp}_{F_i(\Phi)}$ be the exponential map of $U_i$ at $F_i(\phi)$. Note that the exponential map on $T_{\Phi}Z_j$ is the restriction of $F_i^{-1}\circ\mathrm{exp}_{F_i(\Phi)}\circ F_i$ to the subspace $T_{\Phi}Z_j$, for all $i$ and $j$. Thus for $X\in \bigcup_{j} T_{\Phi}Z_j$, we can define $\mathrm{exp}_{\Phi}(X):=F_i^{-1}(\mathrm{exp}_{F_i(\Phi)}(F_i(X)))$, which is independent of the choice of $i\in\{1,2,3\}$.

In practice, given any $X\in T$, we can compute $\mathrm{exp}_{F_i(\Phi)}$ with Equation (\ref{exp}). Then we will check if the following equalities hold: 
\begin{equation}\label{Tgood}
F_1^{-1}(\mathrm{exp}_{F_1(\Phi)}(F_1(X)))=F_2^{-1}(\mathrm{exp}_{F_2(\Phi)}(F_2(X)))=F_3^{-1}(\mathrm{exp}_{F_3(\Phi)}(F_3(X))).
\end{equation}
If (\ref{Tgood}) holds, then $X\in T_{\Phi}Z_j$ for some $j$ and we obtain a perfect tensor $\mathrm{exp}_{\Phi}(X)\in Z_j\subset\PP(4,d)$. Pick some linearly independent vectors $X_k\in T$, and let $t_k$'s be real symbolic variables. If (\ref{Tgood}) holds for $X=\sum_{k} t_kX_k$, then $\mathrm{span}_{\BR}\{X_k\mid k\}\in T_{\Phi}Z_j$ for some $j$, and we obtain a positive-dimensional family $\mathrm{exp}(\mathrm{span}_{\BR}\{X_j\mid j\})\subset Z_j$ parameterized by $t_j$'s.

Since $F_i$ is a linear isomorphism, $T_{\Phi}F_i^{-1}(U_i)=F_i^{-1}(T_{F_i(\Phi)}U_i)$. The equation for $T_{F_i(\Phi)}U_i$ is $XF_i(\Phi)^{\dagger}+F_i(\Phi)X^{\dagger}=0$. Thus, the equations for the space $T$ is
\begin{equation}\label{Teqn}
F_i(X)F_i(\Phi)^{\dagger}+F_i(\Phi)F_i(X)^{\dagger}=0,\forall i=1,2,3.
\end{equation}

To summarize, our method to find families of perfect tensors in $\PP(4,d)$ is:
\begin{itemize}
    \item Step 1: pick a known solution $\Phi\in\PP(4,d)$;
    \item Step 2: compute $T$ via equation (\ref{Teqn});
    \item Step 3: pick some linearly independent vectors $X_k$'s in $T$, then verify if (\ref{Tgood}) holds when for $X=\sum_k t_kX_k$. If it holds, we obtain $\mathrm{exp}_{\Phi}(\mathrm{span}_{\mathbb{R}}\{X_j\mid j\})$, a parameterized family of perfect tensors in $\PP(4,d)$.
\end{itemize}

\begin{rem}
In step 1, we usually choose $\Phi$ to be the tensor corresponding to an OLS($d$) if $d\neq 6$. Since $F_i(\Phi)$ is a sparse matrix for all $i$, it would be much faster to compute $T$ and the exponential maps afterwards.
\end{rem}

\section{Results on $\PP(4,3)$}
We choose $\Phi$ to be the tensor from Example \ref{OLS3}. We compute the intersection of the tangent spaces $T$ with Python and then perform step 3 with MATLAB, and obtain the following results for $\PP(4,3)$.

\begin{prop}\label{resultintersection}
The triple intersection of the tangent spaces $T:=\bigcap_{i=1}^3 T_{\Phi}F_i^{-1}(U_i)$ equals to all 3 pairwise intersections $T_{\Phi}F_i^{-1}(U_i)\cap T_{\Phi}F_j^{-1}(U_j)$, $\forall 1\leq i<j\leq3$.
\end{prop}

\begin{prop}\label{resultbasis} $T$ is a 33-dimensional real vector space. A basis is given below, where $\mathbf{i}$ stands for the imaginary unit, and $ijkl$ denotes the computational basis $|ijkl\rangle$:

{\small
\begin{tabular}{cc}
$e_1=-1132+1223-2113+2231-3121+3212,$ &$f_1=(1132+1223+2113+2231+3121+3212)\mathbf{i},$\\
$e_2=-1111+1323-2122+2331-3133+3312,$ &$f_2=(1111+1323+2122+2331+3133+3312)\mathbf{i},$\\
$e_3=-1211+1332-2222+2313-3233+3321,$ &$f_3=(1211+1332+2222+2313+3233+3321)\mathbf{i},$\\
$e_4=-1131-1213-1322+2123+2232+2311,$ &$f_4=(1131+1213+1322+2123+2232+2311)\mathbf{i},$\\
$e_5=-1112-1221-1333+3123+3232+3311,$ &$f_5=(1112+1221+1333+3123+3232+3311)\mathbf{i},$\\
$e_6=-2112-2221-2333+3131+3213+3322,$ &$f_6=(2112+2221+2333+3131+3213+3322)\mathbf{i},$\\
$e_7=1113-1321-2223+2312-3122+3211, $ &$f_7=(1113+1321+2223+2312+3122+3211)\mathbf{i},$\\
$e_8=-1133+1222+2121-2332-3231+3323,$ &$f_8=(1133+1222+2121+2332+3231+3323)\mathbf{i},$\\
$e_9=1212-1331+2111-2233-3132+3313, $ &$f_9=(1212+1331+2111+2233+3132+3313)\mathbf{i},$\\
$e_{10}=1122-1233+2212-2323-3113+3332,$ &$f_{10}=(1122+1233+2212+2323+3113+3332)\mathbf{i},$\\
$e_{11}=1121-1313-2133+2211-3223+3331,$ &$f_{11}=(1121+1313+2133+2211+3223+3331)\mathbf{i}$,\\
$e_{12}=-1231+1312+2132-2321-3111+3222$, &$f_{12}=(1231+1312+2132+2321+3111+3222)\mathbf{i}$,\\
\multicolumn{2}{c}{$g_1$ to $g_9$: $1123\mathbf{i}, 1232\mathbf{i}, 1311\mathbf{i}, 2131\mathbf{i}, 2213\mathbf{i},2322\mathbf{i}, 3112\mathbf{i},3221\mathbf{i}, 3333\mathbf{i}.$}
\end{tabular}
}%
\end{prop}

This basis consists of 12 pairs of tensors of rank 6, and 9 pure imaginary tensors of rank 1. For each $j$, $e_j$ and $f_j$ form a pair, in the sense that they have the same support (i.e., the set of computational basis vectors $|ijkl\rangle$ that appear in the sum), and $e_j$ is pure real while $f_j$ is pure imaginary. Any two vectors in this basis have disjoint supports, unless they are a pair $(e_j,f_j)$.

\begin{prop}\label{result4dimension}  All vectors in the following 12 4-dimensional spaces satisfy Equation (\ref{Tgood}):
$$\mathrm{span}_{\mathbb{R}}\{e_i,e_j,f_i,f_j\},\text{ for }1\leq i<j\leq 3,4\leq i<j\leq 6,7\leq i<j\leq 9,\text{ or }10\leq i<j\leq 12,$$
which give us 12 4-dimensional families of perfect tensors $\mathrm{exp}_{\Phi}(\mathrm{span}_{\mathbb{R}}\{e_i,e_j,f_i,f_j\})\subset\PP(4,3)$. 
\end{prop}

\begin{ex}\label{ex1}
We obtain the explicit parameterization of the family $\mathrm{exp}_{\Phi}(\mathrm{span}_{\mathbb{R}}\{e_1,e_2,f_1,f_2\})$. Write $\Psi:=\mathrm{exp}_{\Phi}(t_1e_1+t_2f_1+t_3e_2+t_4f_2)=\sum_{i,j,k,l=1}^{3}\Psi_{ijkl}|ijkl\rangle$ and $||t||:=(t_1^2+t_2^2+t_3^2+t_4^2)^{\frac{1}{2}}$, then
\begin{align*}
    \Psi_{1111}&=\Psi_{2122}=\Psi_{3133}=-\mathrm{sinh}(||t||\mathbf{i})(t_3 - \mathbf{i}t_4)/(||t||\mathbf{i}),\\
    \Psi_{1211}&=\Psi_{2222}=\Psi_{3233}=(\mathrm{cos}(||t||)-1)(t_1 + t_2\mathbf{i})(t_3 - t_4\mathbf{i})/||t||^2,\\
    \Psi_{1311}&=\Psi_{2322}=\Psi_{3333}=((t_3^2+t_4^2)\mathrm{cos}(||t||) + t_1^2 + t_2^2)/||t||^2,\\
    \Psi_{3112}&=\Psi_{1123}=\Psi_{2131}=\mathrm{cosh}(||t||\mathbf{i}),\\
    \Psi_{3212}&=\Psi_{1223}=\Psi_{2231}=\mathrm{sinh}(||t||\mathbf{i})(t_1 + \mathbf{i}t_2)/(||t||\mathbf{i}),\\
    \Psi_{3312}&=\Psi_{1323}=\Psi_{2331}=\mathrm{sinh}(||t||\mathbf{i})(t_3 + \mathbf{i}t_4)/(||t||\mathbf{i}),\\
    \Psi_{2113}&=\Psi_{3121}=\Psi_{1132}=-\mathrm{sinh}(||t||\mathbf{i})(t_1 - \mathbf{i}t_2)/(||t||\mathbf{i}),\\
    \Psi_{2213}&=\Psi_{3221}=\Psi_{1232}=((t_1^2+t_2^2)\mathrm{cos}(||t||) + t_3^2 + t_4^2)/||t||^2,\\
    \Psi_{2313}&=\Psi_{3321}=\Psi_{1332}=(\mathrm{cos}(||t||)-1)(t_1 - t_2\mathbf{i})(t_3 + t_4\mathbf{i}))/||t||^2,
\end{align*}
and all $\Psi_{ijkl}$'s not listed above are zero.
\end{ex}

\begin{thm}
For each one of the 12 families of perfect tensors given in Proposition \ref{result4dimension}, there is a 4-dimensional subset consisting of perfect tensors not locally equivalent to any solution corresponding to any OLS(3). 
\end{thm}
\begin{proof}
For any perfect tensor $\Psi$, a tensor is locally equivalent to $\Psi$ if and only if it belongs to the orbit $\mathrm{Orb}_{\Psi}:=(U_1\times U_2\times U_3\times U_4)\cdot \Psi$, where $U_i$ is the unitary group on $\HH_i$. Without loss of generality, we only need to show the existence of a 4-dimensional subset of $\mathrm{exp}_{\Phi}(\mathrm{span}_{\mathbb{R}}\{e_1,e_2,f_1,f_2\})$ that has no intersection with $\mathrm{Orb}_{\Psi}$ for any $\Psi$ corresponding to an OLS(3).

Note that the orbit $\mathrm{Orb}_{\Phi}$ is smooth at $\Phi$, so $\mathrm{Orb}_{\Phi}=\mathrm{exp}_{\Phi}(T_{\Phi}\mathrm{Orb}_{\Phi})$. Let $\mathfrak{u}_i$ be the Lie algebra of $U_i$. With the relation
$$T_{\Phi}\mathrm{Orb}_{\Phi}=T_{\Phi}((U_1\times U_2\times U_3\times U_4)\cdot\Phi)=(\mathfrak{u}_1\cdot\Phi)\oplus\cdots\oplus(\mathfrak{u}_4\cdot\Phi)$$
and the fact that the equation $\sum_{i=1}^4 (X_i\cdot\Phi)=e_1$ has no solutions for $X_i\in \mathfrak{u}_i$, we conclude that $e_1\not\in T_{\Phi}\mathrm{Orb}_{\Phi}$. Therefore any tangent vector $t_1e_1+t_2f_1+t_3e_2+t_4f_2$ with $t_1\neq 0$ does not belong to $T_{\Phi}\mathrm{Orb}_{\Phi}$, which implies the 4-dimensional subset $S:=\mathrm{exp}_{\Phi}(\mathrm{span}_{\mathbb{R}}\{e_1,e_2,f_1,f_2\}\backslash\mathrm{span}_{\mathbb{R}}\{e_2,f_1,f_2\})$ has no intersection with $\mathrm{Orb}_{\Phi}$.

Since there are only finitely many OLS(3)'s, there exists a neighborhood $N$ of $\Phi$ that has no intersection with $\mathrm{Orb}_{\Psi}$ for any $\Psi$ corresponding to an OLS(3) but not equivalent to $\Phi$. We obtain is a 4-dimensional subset $N\cap S\subset\mathrm{exp}_{\Phi}(\mathrm{span}_{\mathbb{R}}\{e_1,e_2,f_1,f_2\})$ consisting of perfect tensors not corresponding to any OLS(3).
\end{proof}

\begin{prop}\label{resultclassical} All vectors in $\mathrm{span}_{\mathbb{R}}\{g_j\mid 1\leq j\leq 9\}$ satisfy Equation (\ref{Tgood}). For the general vectors $X=\sum_{j=1}^9 t_jg_j$ in this subspace, $F_1(\mathrm{exp}_{\Phi}(X))=$
$$\begin{pmatrix}
  0&   0&   0&   0&   0&  \mathrm{exp}(\mathbf{i}t_1)&   0&   0&   0\\
  0&   0&   0&   0&   0&   0&   0& \mathrm{exp}(\mathbf{i}t_2)&   0\\
\mathrm{exp}(\mathbf{i}t_3)&   0&   0&   0&   0&   0&   0&   0&   0\\
  0&   0&   0&   0&   0&   0& \mathrm{exp}(\mathbf{i}t_4)&   0&   0\\
  0&   0& \mathrm{exp}(\mathbf{i}t_5)&   0&   0&   0&   0&   0&   0\\
  0&   0&   0&   0& \mathrm{exp}(\mathbf{i}t_6)&   0&   0&   0&   0\\
  0& \mathrm{exp}(\mathbf{i}t_7)&   0&   0&   0&   0&   0&   0&   0\\
  0&   0&   0& \mathrm{exp}(\mathbf{i}t_8)&   0&   0&   0&   0&   0\\
  0&   0&   0&   0&   0&   0&   0&   0& \mathrm{exp}(\mathbf{i}t_9)\\
\end{pmatrix}$$
which gives us a 9-dimensional family of $\PP(4,3)$. 
\end{prop}

Note that unlike the families given in Proposition \ref{result4dimension}, $\mathrm{exp}_{\Phi}(\sum_{j=1}^9 t_jg_j)$ is locally equivalent to $\Phi$ for all $t_j\in\mathbb{R}$.

\begin{prop}\label{result6dimension}
\begin{enumerate}[{(1)}]
    \item For a general vector in the following four 6-dimensional subspaces, the first 13th terms of the Taylor series of the three exponential maps in Equations (\ref{Tgood}) agree:
    \begin{equation}\label{6dimensionalsubspace}
    \mathrm{span}_{\mathbb{R}}\{e_i,e_{i+1},e_{i+2},f_i,f_{i+1},f_{i+2}\},\text{ for }i=1,4,7,10.   
    \end{equation}
    \item For any vector in $T=\bigcap_{i=1}^3 T_{\Phi}F_i^{-1}(U_i)$ but not in the four 6-dimensional subspaces in (\ref{6dimensionalsubspace}) or the 9-dimensional subspace $\mathrm{span}_{\mathbb{R}}\{g_j\mid 1\leq j\leq 9\}$, the Taylor series of the three exponential maps in Equations (\ref{Tgood}) does not agree from degree 2, so $\mathrm{exp}_{\Phi}$ is not defined at those vectors.
\end{enumerate}

\end{prop}

In the first statement of Proposition \ref{result6dimension}, it is not known if further terms in Taylor series agree not not, because it takes too much time to compute terms of higher degrees on my computer. If all terms agree, i.e., the Equations (\ref{Tgood}) hold for all four 6-dimensional subspaces in (\ref{6dimensionalsubspace}), then there will be two good results:
\begin{enumerate}
    \item we will obtain four 6-dimensional families of perfect tensors $$\mathrm{exp}_{\Phi}(\mathrm{span}_{\mathbb{R}}\{e_i,e_{i+1},e_{i+2},f_i,f_{i+1},f_{i+2}\}),$$ 
    which will be a better result than Proposition \ref{result4dimension};
    \item together with Proposition \ref{resultclassical} and the second statement of Proposition \ref{result6dimension}, the following five subspaces of $T$ will be exactly all maximal subspaces of $T$ where Equation (\ref{Tgood}) holds:
    \begin{equation}
    \mathrm{span}_{\mathbb{R}}\{g_j\mid 1\leq j\leq 9\}, \mathrm{span}_{\mathbb{R}}\{e_i,e_{i+1},e_{i+2},f_i,f_{i+1},f_{i+2}\}), \text{ for }i=1,4,7,10.
    \end{equation}
    Therefore the all $T_{\Phi}Z_j$'s will be the five subspaces above, and all $Z_j$'s will be the exponential of these subspaces.
\end{enumerate}

\begin{conj}
There exist exactly five maximal irreducible subvarieties of $\PP(4,3)$ that are smooth at $\Phi$ that are smooth at $\Phi$. They are:
$$\mathrm{exp}_{\Phi}(\mathrm{span}_{\mathbb{R}}\{g_j\mid 1\leq j\leq 9\}),\mathrm{exp}_{\Phi}(\mathrm{span}_{\mathbb{R}}\{e_i,e_{i+1},e_{i+2},f_i,f_{i+1},f_{i+2}\}),\text{ for }i=1,4,7,10.$$
\end{conj}

\section{Results on $\PP(4,4)$ and $\PP(4,5)$}
Our choice of $\Phi$ is given by the following OLS's. 
\begin{center}
\begin{tabular}{ |c|c|c|c| }
\hline
 1,1 & 2,3 & 3,4 & 4,2\rule[-1ex]{0pt}{4ex}\\
\hline
 2,2 & 1,4 & 4,3 & 3,1\rule[-1ex]{0pt}{4ex}\\
\hline
 3,3 & 4,1 & 1,2 & 2,4\rule[-1ex]{0pt}{4ex}\\
\hline
 4,4 & 3,2 & 2,1 & 1,3\rule[-1ex]{0pt}{4ex}\\
\hline
\end{tabular}\qquad
\begin{tabular}{ |c|c|c|c|c| }
\hline
 1,1 & 2,4 & 3,2 & 4,5 & 5,3\rule[-1ex]{0pt}{4ex}\\ 
\hline
 2,2 & 3,5 & 4,3 & 5,1 & 1,4\rule[-1ex]{0pt}{4ex}\\
\hline
 3,3 & 4,1 & 5,4 & 1,2 & 2,5\rule[-1ex]{0pt}{4ex}\\
\hline
 4,4 & 5,2 & 1,5 & 2,3 & 3,1\rule[-1ex]{0pt}{4ex}\\
\hline
 5,5 & 1,3 & 2,1 & 3,4 & 4,2\rule[-1ex]{0pt}{4ex}\\
\hline
\end{tabular}
\end{center}

With the same method to study $\PP(4,3)$, we obtain the following results.

\begin{prop} The statement of Proposition \ref{resultintersection} holds for $\PP(4,5)$ but not for $\PP(4,4)$.
\end{prop}

\begin{prop normal} For $\PP(4,5)$, $T$ is 145-dimensional. Similar to Proposition \ref{resultbasis}, a basis of $T$ consists of 60 pairs of vectors $(e_j,f_j)$ of rank 10 and 25 pure imaginary vectors $g_k$ of rank 1. 

For $\PP(4,4)$, $T$ is 76-dimensional. Unlike $\PP(4,3)$ and $\PP(4,5)$, a basis of $\PP(4,4)$ consists of 24 pairs vectors $(e_j,f_j)$'s of rank 8, 16 pure imaginary vectors $g_k$'s of rank 1, and 12 pure imaginary tensors $h_l$'s of rank 4. 

In either $\PP(4,4)$ or $\PP(4,5)$, these basis vectors have disjoint supports except the pairs $(e_j,f_j)$. While the Equation (\ref{Tgood}) holds for all $e_j$'s, $f_j$'s and $g_k$'s in $\PP(4,4)$ and $\PP(4,5)$, it does not hold for any $h_l$'s in $\PP(4,4)$.
\end{prop normal}

\begin{prop} For $\PP(4,4)$ and $\PP(4,5)$, all vectors in $\mathrm{span}_{\mathbb{R}}\{g_j\mid j\}$ satisfy Equation (\ref{Tgood}). We obtain a 16-dimensional family in $\PP(4,4)$ and a 25-dimensional family in $\PP(4,5)$ of the form $\mathrm{exp}(\sum_j t_jg_j)$. Perfect tensors in these two families are locally equivalent to the two tensors corresponding to the two OLS's given at the beginning of this section respectively.
\end{prop}

\section*{Acknowledgements}
I appreciate my advisor Joseph Landsberg for instructions on my research on perfect tensors, and comments and corrections to this paper. I also thank Ion Nechita for useful conversations and Wojciech Bruzda for suggestions on simplifying the expressions in Example \ref{ex1}.

\bibliographystyle{amsplain}
\bibliography{perfect.bib}
\end{document}